\documentclass{amsart}

\usepackage{amsfonts}
\usepackage{cases}
\usepackage{mathrsfs}
\usepackage{bbm}
\usepackage{amssymb}
\usepackage{txfonts}
\usepackage{amscd}
\usepackage{amsfonts,latexsym,amsmath,amsthm,amsxtra,mathdots}
\usepackage[bookmarks,colorlinks]{hyperref}
\usepackage[all,cmtip]{xy}
\RequirePackage{amsmath} \RequirePackage{amssymb}
\usepackage{color}
\usepackage{colordvi}
\usepackage{multicol}
\usepackage{pdfsync}
\usepackage[utf8]{inputenc}
\usepackage{hyperref}
\usepackage{graphicx}
\usepackage{amsmath}
\usepackage{amsmath,amscd}
\usepackage[normalem]{ulem}
\usepackage{textcomp}
\usepackage{enumitem}
\usepackage{tikz}
\usetikzlibrary{matrix,arrows,trees}
\usetikzlibrary{positioning}
\usepackage{pdfcolmk}
\usepackage{todonotes}

    \newcommand{\BE}{{\mathbb {E}}} \newcommand{\BF}{{\mathbb {F}}}

     \newcommand{\BZ}{{\mathbb {Z}}}

     \newcommand{\CF}{{\mathcal {F}}}

\def\-{^{-1}}

\newcommand{\delete}[1]{}

    \theoremstyle{plain}

\newtheorem{thm}{Theorem}

\newtheorem{rem}[thm]{Remark}

    \numberwithin{equation}{section}

\def\Proof{\noindent{\bf Proof}\quad}
\def\qed{\hfill$\square$\smallskip}

\begin{document}

\title{Note on the injectivity of the Loday assembly map}


\author{Mark Ullmann}
\address{Freie Universit\"at Berlin, Institut f\"ur Mathematik, Arnimallee 7, 14195 Berlin, Germany}
\email{mark.ullmann@math.fu-berlin.de}

\author{Xiaolei Wu}
\address{Max Planck Institute for Mathematics, Vivatsgasse 7, 53111 Bonn, Germany }
\email{hsiaolei.wu@mpim-bonn.mpg.de}





\begin{abstract}
We show that for finite groups the Loday assembly map with coefficients in finite fields is in general not injective.
\end{abstract}

\maketitle

Given a group $G$, a ring $R$, one can consider the following Loday assembly map \cite[Proposition 4.1.1]{Lo}:

$$\alpha: H_{\ast} (BG,K(R))  \longrightarrow K_{\ast}(R[G])$$

In general, the map is far from isomorphism. For example, the nonvanishing of whitehead torsion will imply the assembly is  not surjective. But when the group is torsion free and the ring is regular, it is implied by the Farrell--Jones Conjecture that this is an isomorphism, see \cite{Lueck-Reich-2005} for definitions, and in particular Corollary 67 for the isomorphism result.  One might ask, would the map be injective for groups with torsion and regular rings, see for example  \cite[Question 1.1]{Ma}. We show this is not the case even if the group is finite  and $R$ is a field.

To construct the counterexample, we  need Quillen's calculation for algebraic K-theory of finite fields.
\begin{thm}[{\cite[Theorem 8]{Qu}}]\label{Qu}
Let $\BF_q$ be a finite field with $q$ elements, then for $n >0$, $K_n(\BF_q) = 0$ when $n$ is even, and  $K_n(\BF_q)$ is a cyclic group of order $q^i -1$ when $n =2i-1$.
\end{thm}
 Note that since $\BF_q$ is a field, $K_n(\BF_q) = 0$ when $n < 0$ and $K_0(\BF_q) = \BZ$. 

\begin{thm}
Let $G$ be a finite group such that $H_2(G,\BZ)$ is nontrivial and $\BF$ a finite field with characteristic $p$  which does not divide the order of $G$, then the assembly map $\alpha$ is not injective. For example, we can take $G=\BZ/2\times \BZ/2$ and $\BF$ to be any finite field with characteristic $p>2$.
\end{thm}

\Proof Since the characteristic of $\BF$ does not divide the order of $G$, by Maschke's theorem, $\BF G$ is semisimple. Since it is finite, by Artin--Wedderburn's Theorem, $\BF G$ is a direct product of matrix rings over finite fields of characteristic $p$. Note that for any matrix ring $M_n(\BE)$ over a field $\BE$, $K_\ast(M_n(\BE)) = K_\ast(\BE)$. Thus by Theorem \ref{Qu}, we in particular have $K_n(\BF G) = 0$ when $n$ is even and $K_0(\BF G) = \BZ^{d}$, where $d$ is the number of components in the Artin--Wedderburn Theorem for $\BF G$. Note also that $K_n(\BF G) = 0$ when $n<0$.

Now applying Atiyah--Hirzebruch spectral sequence to the left hand side of the assembly map $H_{\ast} (BG,K(\BF))$, $E_{p,q}^2 = H_p(BG,K_q(\BF))$. (abbreviate $H_\ast(G):= H_\ast(G,\BZ)$)

\begin{tikzpicture}
  \matrix (m) [matrix of math nodes,
    nodes in empty cells,nodes={minimum width=5ex,
    minimum height=5ex,outer sep=-5pt},
    column sep=1ex,row sep=1ex]{
           {\scriptstyle q}    &      &     &     & & \\           
           3   &  K_3(\BF) &  \ast  & \ast & \ast& \\
          2    &  0 &  0  & 0 & 0& \\
          1     & K_1(\BF) & \ast  & \ast & \ast& \\
          0     & K_0(\BF)  & H_1(G) &  H_2(G)  & H_3(G) & \\
    \quad\strut &   0  &  1  &  2  & 3 & \strut & {\scriptstyle p} \\};
\draw[thick] (m-1-1.east) -- (m-6-1.east) ;
\draw[thick] (m-6-1.north) -- (m-6-6.north) ;
\end{tikzpicture}

Note that the assembly map for $K_0$ and $K_1$ is indeed injective and the terms $E_{0,0}^2, E_{1,0}^2, E_{0,1}^2$ survive \cite[Lemma 2]{Lueck-Reich-2005}. Thus $H_2(BG,K_0(\BF)) = H_2(G)$ survives in the spectral sequence. But since $H_2(G)\not = 0$ and $K_2(\BF G) = 0$, the assembly map can not be injective. \qed

\begin{rem}
It would be interesting if one could find groups $G$ such that the assembly map is not injective for the ring $\BZ$. However it seems that there is not much calculation known in these cases. 
\end{rem}

\begin{rem}
 In general, the Farrell--Jones Conjecture predicts that the assembly map
 \begin{equation*}
    H^G_q(E_{\CF} G, K_{\mathbf{R}} ) \to K_q(R[G])
 \end{equation*}
is an isomorphism for any group $G$, any ring $R$ and $\CF$ the family of virtually cyclic subgroups of $G$ \cite{Lueck-Reich-2005}.
 There are results that the family can be made smaller, see \cite[Section 7]{Bartels-Farrell-Lueck-cocompact} and  \cite{DQR}.   The Transitivity Principle \cite[Theorem 65]{Lueck-Reich-2005} and the inheritance under subgroups implies, that if the Farrell--Jones Conjecture holds for a family $\CF$, it holds for all families containing $\CF$.  There are also injectivity results for smaller families, for example for the family of finite groups, cf.~\cite{Kasprowski-FDC}.  Thus one may think, that the relative assembly map to a larger family always gives an injective map.  Our example shows, that this is not true for the trivial family to the family of all (or finite) subgroups.
 \end{rem}

\textbf{Acknowledgements.} The first author was supported by the DFG--project C1 of CRC 647 in Berlin. The second author is supported by the Max Planck Institute for Mathematics at Bonn.


\begin{thebibliography}{1}
{\small 

 \bibitem{Bartels-Farrell-Lueck-cocompact}
 A.~Bartels, F.~T. Farrell, and W.~L{\"u}ck.
 \newblock The {F}arrell-{J}ones conjecture for cocompact lattices in virtually
   connected {L}ie groups.
\newblock {\em J. Amer. Math. Soc.}, 27(2):339--388, 2014.

 \bibitem{DQR} J. F. Davis, F. Quinn and H. Reich. Algebraic K-theory over the infinite dihedral group: a controlled topology approach, J. Topology,   4 (2011) 505-528.
 
\bibitem{Kasprowski-FDC}
D. Kasprowski.
\newblock On the {$K$}-theory of groups with finite decomposition complexity.
\newblock {\em Proc. Lond. Math. Soc. (3)}, 110(3):565--592, 2015.


\bibitem{Lueck-Reich-2005}
W.~L{\"u}ck and H.~Reich.
\newblock The {B}aum-{C}onnes and the {F}arrell-{J}ones conjectures in {$K$}-
  and {$L$}-theory.
\newblock In {\em Handbook of $K$-theory. Vol. 1, 2}, pages 703--842. Springer,
  Berlin, 2005.
  
\bibitem{Lo} J.-L. Loday, K-th\'eorie alg\'ebrique et repr\'esentations de groupes, Ann. Sci. \'Ecole Norm. Sup. (4) 9 (1976), no. 3, 309-377.  

\bibitem{Ma}  C. Malkiewich.
\newblock Coassembly and the K-theory of finite groups. 
\newblock arXiv:1503.06504


\bibitem{Qu} D. Quillen.
\newblock On the cohomology and K-theory of the general linear groups over a finite field.  
\newblock {\em Ann. of Math. (2) 96} (1972), 552-586. 



}
\end{thebibliography}
\end{document}